\documentclass[english,12pt]{amsart}
\pdfoutput=1
\usepackage[T1]{fontenc}
\usepackage[latin9]{inputenc}
\usepackage[letterpaper]{geometry}
\usepackage{babel}

\usepackage{verbatim}
\usepackage{amsthm}
\usepackage{graphicx}
\usepackage[unicode=true, pdfusetitle,
 bookmarks=true,bookmarksnumbered=false,bookmarksopen=false,
 breaklinks=false,pdfborder={0 0 0},backref=false,colorlinks=false]
 {hyperref}

\makeatletter
\numberwithin{equation}{section} 
\numberwithin{figure}{section} 
\theoremstyle{plain}
\theoremstyle{plain}
\newtheorem{thm}{Theorem}
  \theoremstyle{plain}
  \newtheorem{prop}[thm]{Proposition}
  \theoremstyle{definition}
  \newtheorem{defn}[thm]{Definition}
  \theoremstyle{plain}
  \newtheorem{cor}[thm]{Corollary}
  \theoremstyle{plain}
  \newtheorem{lem}[thm]{Lemma}

\makeatother

\begin{document}

\title{A synthesis for exactly 3-edge-connected graphs}

\author{Carl Kingsford}

\address{C. Kingsford, Department of Computer Science and Institute for Advanced
Computer Studies, University of Maryland, College Park, MD, USA}

\email{carlk@cs.umd.edu}

\thanks{C.K. was partially supported by NSF grant IIS-0812111.}

\author{Guillaume Mar\c{c}ais}

\address{G. Mar\c{c}ais, Program in Applied Mathematics \& Statistics and Scientific
Computation, University of Maryland, College Park, MD, USA}

\email{guillaume@marcais.net}

\date{\today}

\begin{abstract}
A multigraph is \textit{exactly k-edge-connected} if there are exactly
$k$ edge-disjoint paths between any pair of vertices. We characterize
the class of exactly 3-edge-connected graphs, giving a synthesis involving
two operations by which every exactly 3-edge-connected multigraph
can be generated. Slightly modified syntheses give the planar exactly
3-edge-connected graphs and the exactly 3-edge-connected graphs with
the fewest possible edges.
\end{abstract}
\maketitle

\section{Introduction}

We define a multigraph $G$ to be \textbf{exactly} $k$-edge-connected
if there are exactly $k$ edge-disjoint paths between any pair of
distinct vertices $u,v$ in $G$. The study of exactly $k$-edge-connected
graphs has several motivations. When designing a communication network,
a natural requirement is that every pair of vertices can communicate
over $\geq k$ edge-disjoint pathways. This ensures the network is
robust to attacks and edge failures. Exactly connected networks are
the most fair because no pair of vertices is provided additional communication
pathways. In addition, if two antagonistic agents arise in an exactly
$k$-edge-connected network of formerly mutually friendly agents,
communication between the two antagonistic agents can always be severed
by interrupting $k$ edges, while still guaranteeing that, before
the split, communication ability was robust. 

Exactly connected graphs arise in other contexts as well. For example,
they are obtained if all pairs of vertices of higher local edge-connectivity
are merged. Suppose $G$ is a $k$-edge-connected graph, but not necessarily
exactly $k$-edge-connected. If the sets of vertices that are mutually
connected by $>k$ edge-disjoint paths are collapsed into supernodes,
then the resulting multigraph $G'$ is either a graph of a single
vertex or is exactly $k$-edge-connected. This follows because $G'$
is at least $k$-edge-connected as $G$ was, and no two supernodes
can be connected by $>k$ edge-disjoint paths, otherwise they would
have been merged. That the supernodes are equivalence classes follows
directly from Menger's theorem~\cite{diestel05}.

\begin{figure}
\begin{centering}
\includegraphics[clip,width=0.3\textwidth]{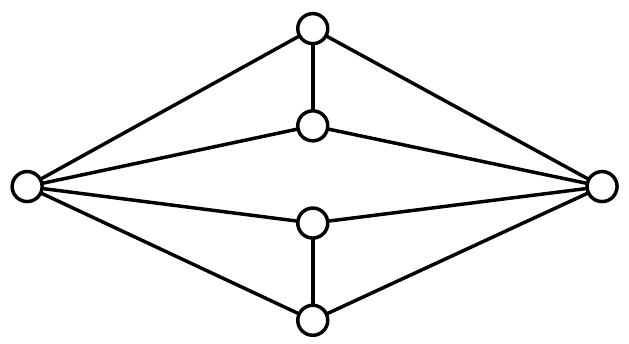}
\par\end{centering}

\caption{\label{fig:lips}A 3-edge-connected graph that is edge minimal but
not exactly 3-edge connected.}

\end{figure}

All exactly $k$-edge-connected graphs are both edge-minimal and edge-maximal
in the sense that removing or adding an edge will destroy exact $k$-edge
connectivity. In the same sense, they are vertex-minimal and vertex-maximal.
On the other hand, not all edge-minimal $k$-edge-connected graphs
are exactly $k$-edge-connected (Figure~\ref{fig:lips}). Hence,
the class of exactly $k$-edge-connected graphs is a strict subset
of the edge-minimal $k$-edge-connected graphs, and their investigation
is justified in its own right. The Harary graphs~\cite{Gross2006},
for example, are exactly $k$-edge-connected, but there is only one
Harary graph for each pair $(n,k)$, where $n$ is the order of the
graph and $k$ the connectivity requirement.

The exactly 2-edge-connected graphs are {}``trees of cycles.'' More
formally, if $G$ is exactly 2-edge-connected, it is the union of
simple cycles such that the block-cutvertex graph $H=(V,E)$ is a
tree, where $V$ contains a vertex $u_{c}$ for each cycle $c$ and
a vertex $u_{v}$ for each vertex $v\in G$ shared by more that one
cycle, and where $E=\{\{u_{c},u_{v}\}:v\in c\}$. Every 2-edge-connected
graph can be generated from a Whitney-Robbins synthesis~\cite{Gross2006}.
The exactly 2-connected graphs can be generated via a {}``Robbins
synthesis'':
\begin{prop}
An exactly 2-edge-connected graph $G$ is obtained from a Robbins
synthesis. In other words, there exists a sequence of graph $G_{i},0\le i\le l$
such that $G_{0}$ is a simple cycle, $G_{l}=G$ and $G_{i+1}$ is
obtained from $G_{i}$ by a cycle addition.\end{prop}
\begin{proof}
$G$ is 2-edge-connected, so it is the result of a Whitney-Robbins
synthesis. Let $G_{i}$ $(0\le i\le l)$ such that $G_{0}$ is a simple
cycle and $G_{i+1}$ is obtained from $G_{i}$ by a path or a cycle
addition. Note that $G_{0}$ is exactly 2-edge-connected. Suppose
that $G_{i}$ is exactly 2-edge-connected and that $G_{i+1}$ is obtained
by a path addition, between $u$ and $v$. There are 2 edge-disjoint
paths between $u$ and $v$ in $G_{i}$. The newly added path is obviously
edge-disjoint from the paths in $G_{i}$. So there are 3 edge-disjoint
paths in $G_{i+1}$ between $u$ and $v$. Given that no edge or vertex
is ever removed, there are also 3 edge-disjoint paths in $G$ between
$u$ and $v$, which is a contradiction. On the other hand, cycle
addition does not create any new edge-disjoint paths between existing
vertices. So $G$ is obtained from $G_{0}$ by a sequence of cycle
additions.
\end{proof}
The situtation is more complicated for $k=3$. Tutte~\cite{tutte66}
gave a synthesis to generate all 3-vertex-connected graphs starting
from a $n$-spoke wheel $W_{n}$ via a series of arbitrary edge additions
and an operation that converts vertices of degree $\geq4$ into edges.
Unfortunately, neither of Tutte's synthesis operations preserve exact
edge-connectivity, and so the question of whether such a synthesis
can be found for exactly $3$-edge-connected graphs is an interesting
one. A computational search reveals that, in fact, there are many
exactly $3$-edge-connected graphs. For example, there are $717$
simple, biconnected, exactly $3$-edge-connected isomorphism classes
on $11$ vertices. The main result of this paper is a synthesis for
exactly $3$-edge-connected graphs, which we give below in Section~\ref{sec:A-synthesis-for}.
Following the proof of the synthesis, we describe some of the properties
of exact $3$-edge-connected graphs.

\section{A synthesis for exactly 3-edge connected graphs\label{sec:A-synthesis-for}}

We begin by describing several operations that preserve exact 3-edge-connectedness
when applied to any exactly 3-edge-connected graph. Eventually, the
synthesis will exploit only two of these operations: block gluing
and cycle expansion. The other operations will play a role in the
proof of the synthesis and are also interesting in their own right.
In the following, biconnected graphs are graphs which are 2-vertex-connected.
A vertex whose removal would disconnect a connected but not biconnected
graph is called an articulation point. A block is a maximal connected
subgraph that is biconnected. Throughout, all sets should be considered
multisets and all graphs are multigraphs and loopless. The notation
$(u,v)^{r}$ denotes an undirected edge between $u$ and $v$ of multiplicity
$r$ (when $r=1$ it is omitted).

\subsection{Gluing operations}

\begin{figure}
\begin{centering}
\includegraphics[clip,width=0.6\textwidth]{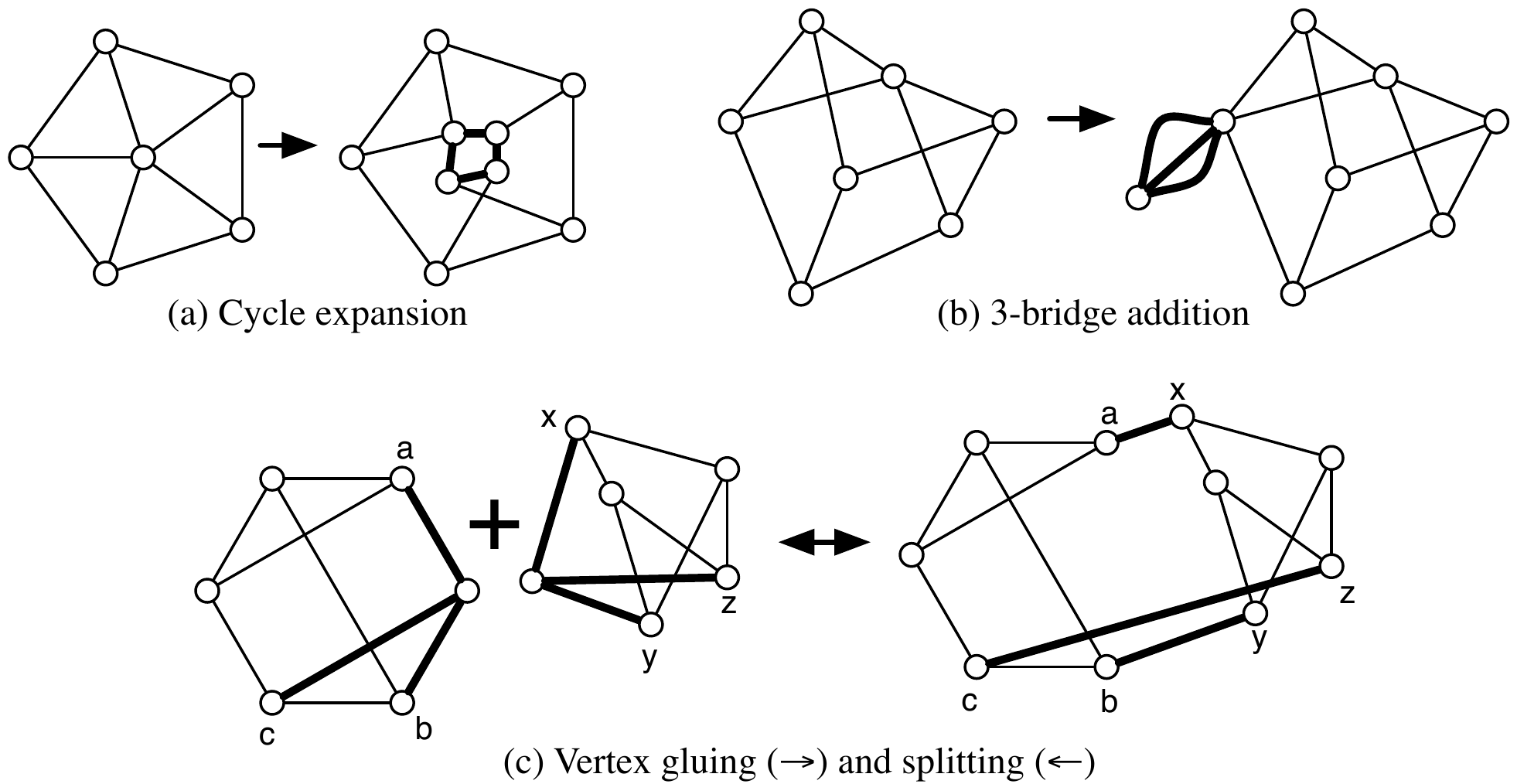}
\par\end{centering}

\caption{\label{fig:Several-operations-that}Several operations that preserve
exact 3-edge-connectivity.}

\end{figure}

In  this section, $G_{1}=(V_{1},E_{1})$ and $G_{2}=(V_{2},E_{2})$
are two exactly $k$-edge-connected graphs with $V_{1}\cap V_{2}=\emptyset$.
We will define operations to glue these graphs together into $G=(V,E)$
that is also exactly $k$-edge-connected. See Figure~\ref{fig:Several-operations-that}
for illustrations of several of the operations. The first collection
of these operations exploits the relationship between exact connectivity
and the blocks, or biconnected components, of a graph.

\begin{defn}
[Block gluing]$G=(V,E)$ is the graph obtained by gluing $G_{1}=(V_{1},E_{1})$
and $G_{2}=(V_{2},E_{2})$ at one vertex. Formally, let $u_{1}\in V_{1}$
and $u_{2}\in V_{2}$, and let $X=\{(u_{1},v):v\in N(u_{1})\}\cup\{(u_{2},v):v\in N(u_{2})\}$
be the adjacent edges of $u_{1}$ and $u_{2}$. Then construct $G=(V,E)$
by setting\begin{align*}
V & =V\cup V_{2}\cup\{u\}\setminus\{u_{1},u_{2}\}\\
E & =E_{1}\cup E_{2}\cup\{(u,v):v\in N(u_{1})\cup N(u_{2})\}\setminus X.\end{align*}
This is called \textbf{block gluing} because the vertex $u$ becomes
an articulation point in $G$.\end{defn}
\begin{prop}
If $G_{1}=(V_{1},E_{1})$ and $G_{2}=(V_{2},E_{2})$ are exactly $k$-edge-connected,
then any block gluing $G$ of them is exactly $k$-edge-connected.\end{prop}
\begin{proof}
Any pair of vertices that are in the same block of $G$ are connected
by exactly $k$ edge-disjoint paths because they were before the gluing
and no additional edge-disjoint path can be created by leaving the
block and re-entering it. Let $u$ be the articulation point created
by the block gluing, and let $w\in V_{1},v\in V_{2}$ be a pair of
vertices in different blocks of $G$. Vertices $w,u$ are connected
by $k$ edge-disjoint paths and $u,v$ are also connected by $k$
edge-disjoint paths. Hence, there are $\geq k$ edge-disjoint paths
between $w,v$. Any such path must pass through $u$, so there cannot
be more than $k$ without creating $>k$ edge-disjoint paths between
$w$ and $u_{1}$ in $G_{1}$, which cannot happen. \end{proof}
\begin{cor}
\label{cor:Block-splitting}The subgraph induced by a block of an
exactly $k$-edge-connected graph is exactly $k$-edge-connected.\end{cor}
\begin{defn}
[k-bridge addition] $G=(V,E)$ is obtained from $G_{1}$ by adding
a vertex and $k$ parallel edges from the vertex to an existing vertex
in $G_{1}$. Formally, $V=V_{1}\cup\{u\}$ and $E=E_{1}\cup\{(v,u)^{k}\}$.\end{defn}
\begin{cor}
If $G$ is exactly $k$-edge-connected, then any $k$-bridge addition
preserves exact $k$-edge-connectivity.
\end{cor}
In addition to the above block-based gluings, we have a vertex gluing
that allows two blocks to be merged into one, and the converse operation,
vertex splitting, which breaks a graph at a minimum cut.
\begin{defn}
[Vertex gluing]Let $u_{1}$ and $u_{2}$ be degree $k$ vertices
of $G_{1}$ and $G_{2}$ respectively. Construct $G$ by removing
$u_{1}$ and $u_{2}$ and pairing together the $2k$ edges from $G_{1}$
and $G_{2}$. Formally, let $v_{1}^{i}$ and $v_{2}^{i}$, $1\le i\le k$,
be respectively the neighbors (not necessarily distinct) of $u_{1}$
in $G_{1}$ and $u_{2}$ in $G_{2}$. We then construct $G$ by setting\begin{align*}
V & =V_{1}\cup V_{2}\setminus\{u_{1},u_{2}\}\\
E & =E_{1}\cup E_{2}\cup\left\{ \left(v_{1}^{i},v_{2}^{i}\right):1\leq i\leq k\right\} \setminus\left\{ \left(u_{1},v_{1}^{i}\right),\left(u_{2},v_{2}^{i}\right):1\leq i\leq k\right\} .\end{align*}
The graph $G$ is a \textbf{vertex gluing} of $G_{1}$ and $G_{2}$. \end{defn}
\begin{prop}
If $G_{1}=\left(V_{1},E_{1}\right)$ and $G_{2}=\left(V_{2},E_{2}\right)$
are exactly $k$-edge-connected, then the vertex gluing $G=\left(V,E\right)$
of $G_{1}$ and $G_{2}$ is exactly $k$-edge-connected.\end{prop}
\begin{proof}
Let $u_{1}\in V_{1}$ and $u_{2}\in V_{2}$ be the vertices used in
the vertex gluing. Consider any pair $v,w$ of vertices in G. If $v\in V_{1}$
and $w\in V_{2}$, then there are no more than $k$ edge-disjoint
paths between them because the cut $\left(V_{1}\setminus\left\{ u_{1}\right\} ,V_{2}\backslash\left\{ u_{2}\right\} \right)$
has cardinality $k$. On the other hand, $k$ edge-disjoint paths
exist between them because there were $k$ edge-disjoint paths from
$v$ to $u_{1}$ and from $u_{2}$ to $w$ which can be combined to
create $k$ edge-disjoint paths between $v$ and $w$.

Suppose, instead, $v,w$ are both in $V_{1}$ (the case for $V_{2}$
is symmetric). Then there are $k$ edge-disjoint $v-w$ paths, say
$P_{i},1\le i\le k$, in $G_{1}$. At most $j\le\left\lfloor k/2\right\rfloor $
of these paths passe through vertex $u_{1}$, say $P_{i},1\le i\le j$.
Let $v_{i}^{1}$ and $w_{i}^{1}$, $1\le i\le j$ be the vertices
that respectively precede and follow $u_{1}$ on path $P_{i}$. Let
$v_{i}^{2}$ and $w_{i}^{2}$ be respectively the neighbors of $v_{i}^{1}$
and $w_{i}^{1}$ in $V_{2}$ after vertex gluing. Let $x$ be a vertex
of $G_{2}$ distinct from $u_{2}$. Given that $G_{2}$ is $k$-edge-connected,
there exists $k$ edge disjoint $x-u_{2}$ paths, say $Q_{i},1\le i\le k$,
in $G_{2}$. Because $u_{2}$ is adjacent to every $v_{i}^{2}$ ans
$w_{i}^{2}$, the $Q_{i}$ paths create edge disjoint $x-v_{i}$ and
$x-w_{i}$ paths, which avoid $u_{2}$, for $1\le i\le j$. In turn,
let combine paths $x-v_{i}$ and $x-w_{i}$ to get $j$ edge disjoint
$v_{i}-w_{i}$ paths, $1\le i\le j$, which avoid $u_{2}$ (named
$R_{i}$). By combining the $P_{i}$ and $R_{i}$ paths, for $1\le i\le j$,
we constructed $v-w$ paths in $G$. 

Conversely, there cannot be more than $\lfloor k/2\rfloor$ edge disjoint
$v-w$ paths in $G$ detouring into $G_{2}$ because of the cut $\left(V_{1}\setminus\left\{ u_{1}\right\} ,V_{2}\setminus\left\{ u_{2}\right\} \right)$.
If there were $>k$ edge disjoint $v-w$ paths in $G$, there would
be $>k$ edge dijsoint $v-w$ paths in $G_{1}$. So $G$ is exactly
$k$-edge-connected.
\end{proof}
{}
\begin{defn}
An edge cut $S$ of a graph $G$ is called \textbf{trivial} if one
of the components of $G\setminus S$ is the trivial graph.
\end{defn}

\begin{defn}
[Vertex splitting]Let $G=(V,E)$ be an exactly $k$-edge connected
graph and $S=\langle V_{1},V_{2}\rangle$ be a non-trivial minimum
cut. Construct $G_{1}=(V_{_{1}}\cup\left\{ x_{1}\right\} ,E_{1})$
and $G_{2}=(V_{2}\cup\left\{ x_{2}\right\} ,E_{2})$ by adding two
new vertices $x_{1}$ and $x_{2}$ attached respectively to $G_{1}$
and $G_{2}$ by $k$ new edges to the vertices adjacent to $S$. Formally,
let $S=\left\{ (u_{i},v_{i})\in V_{1}\times V_{2}:1\le i\le k\right\} $
and\begin{eqnarray*}
E_{1} & = & \left(V_{1}\times V_{1}\cap E\right)\cup\left\{ (x_{1},u_{i}):1\le i\le k\right\} \\
E_{2} & = & \left(V_{2}\times V_{2}\cap E\right)\cup\left\{ (x_{2},v_{i}):1\le i\le k\right\} \,.\end{eqnarray*}
The pair $G_{1},G_{2}$ is called a \textbf{vertex splitting} of $G$
with respect to $S$.\end{defn}
\begin{cor}
\label{cor:Vertex-splitting}Let $G$ be an exactly k-edge-connected
graph, $S$ be a non-trivial minimum cut. $G_{1}$ and $G_{2}$ obtained
by vertex splitting of $G$ with respect to $S$ are exactly $k$-edge-connected.
\end{cor}

\subsection{Cycle Contraction and Expansion}

We describe now two additional operations that create and remove cycles
within exactly $k$-edge-connected graphs. The first, cycle expansion,
is the main non-trivial operation of the synthesis. By convention,
in a multigraph $G$ we view a double edge $\left(u,v\right)^{2}$
as a chordless cycle of length 2. If $C=\left(u_{1},\ldots,u_{n},u_{1}\right)$
is a cycle of length $n\ge3$ and there is a double edge $\left(u_{i},u_{i+1}\right)^{2}$,
then the cycle $C$ has a chord between $u_{i}$ and $u_{i+1}$. 
\begin{prop}
[Cycle expansion]\label{lem:Cycle-expansion}Let $G=(V,E)$ be an
exactly 3-edge-connected and biconnected graph and let $u$ be any
vertex of G. Suppose that $u$ has degree $d$. Let $G'$ be created
from $G$ by replacing vertex $u$ with a cycle $C$ of no more than
$d$ vertices where all the vertices of $C$ are connected to at least
one neighbor of $u$, and all but one vertex of $C$ is connected
to exactly one neighbor of $u$. Formally, let $2\le d'\le d$, $v_{1},\ldots,v_{d}$
be the neighbors of $u$ (not necessarily all distinct) and $G'=(V',E')$
with\begin{align*}
V'= & V\cup\{u_{1}\ldots,u_{d'}\}\setminus\{u\}\\
E'= & E\cup\{(u_{i,}v_{i}):1\leq i\leq d'-1\}\cup\{(u_{d'},v_{i}):d'\leq i\leq d\}\setminus\{(u,v_{i}):1\leq i\leq d\}\,.\end{align*}
 Then $G'$ is exactly 3-edge connected.\end{prop}
\begin{proof}
We shall prove, in order, that there are exactly 3 edge-disjoint paths
in $G'$ between these pairs of vertices:
\begin{enumerate}
\item $(x,y)$ in $V\setminus\{u\}$
\item $(x,u_{i})$ for $x\in V\setminus\{u\}$ and $1\le i\le d'$
\item $(u_{i},u_{j})$ for $1\le i<j\le d'$
\end{enumerate}
For the first case, replacing vertex $u$ by a subgraph, does not
create new paths betwen $x$ and $y$, so there are at most 3 edge-disjoint
paths. Because $u$ is not an articulation point, there exists triplets
of edge-disjoint $x-y$ paths where at most 2 go through $u$. These
2 edge-disjoint $x-y$ paths in $G$ can be changed into 2 edge-disjoint
$x-y$ paths in $G'$ by using some of the edges in the expanded cycle.
So there are 3 edge-disjoint $x-y$ paths in $G'$.

For the second case, fix $i\in[1,d]$ and let $P_{1}$, $P_{2}$ and
$P_{3}$ be 3 edge-disjoint $x-u$ paths in $G$. These paths can
be changed to three paths $P_{1}'$, $P_{2}'$, $P_{3}'$ from $x$
to vertices on the expanded cycle, say $u_{j_{1}}$, $u_{j_{2}}$
and $u_{j_{3}}$. If $i$ is equal to $j_{1}$, then we have 3 edge-disjoint
$x-u_{i}$ paths: $P_{1}'$, $P_{2}'$ plus $u_{j_{2}}-u_{j_{1}}$
on the expanded cycle and $P_{3}'$ plus $u_{j_{3}}-u_{j_{1}}$on
the expanded cycle. If $i$ is not equal to $j_{1}$, $j_{2}$ or
$j_{3}$, let consider a $x-v_{j_{1}}-u_{j_{1}}$ path $P$ (which
can be constructed from an $x-v_{j_{1}}-u_{j_{1}}$ cycle in $G$,
and $G$ is 3 edge connected). This path $P$ is distinct from $P_{1}'$,
$P_{2}'$ and $P_{3}'$ as it contains the edge $(v_{j_{1}},u_{j_{1}})$.
If $P$ is edge-disjoint with at least two of $P_{1}'$, $P_{2}'$
and $P_{3}'$, then we are done. Otherwise, let $e$ be the edge closest
to $u_{j_{1}}$ which is common between $P$ and, say, $P_{1}'$.
Create the new path $P'$ equal to $P_{1}'$ from $x$ to $e$, and
equal to $P$ from $e$ to $u_{j_{1}}$. By construction, this is
a $x-u_{j_{1}}$ path which is edge-disjoint with $P_{2}'$ and $P_{3}'$.

Finally, for the last case. For any pair $(i,j)$, there is a cycle
$C_{ij}$ in $G$ containing the vertices $u$, $v_{i}$ and $v_{j}$.
A piece of this cycle can be changed into an $v_{i}-v_{j}$ path $P$
which is edge-disjoint with the expanded cycle. So there are 3 edge-disjoint
$v_{i}-v_{j}$ paths in $G'$: two paths in the expanded cycle and
the path $P$. There cannot be more than 3 edge-disjoint paths as,
by construction, at least one of $u_{i}$ or $u_{j}$ has degree 3.
\end{proof}
If $G$ is not assumed to be biconnected and $u$ is an articulation
point, a block-respecting cycle expansion can be performed at $u$.
A block-respecting cycle expansion is a cycle expansion done within
one block. Formally, these are defined by breaking the graph $G$
into its blocks with block splitting operations, performing a cycle
expansion on $u$ in one of its blocks, and then gluing all the blocks
back together with block gluing operations. Block-respecting cycle
expansion also preserves exact 3-edge-connectivity. They are not needed
for the synthesis or its proof.

See Figure~\ref{fig:Several-operations-that} for an example of cycle
expansion. The following proposition, where a cycle is contracted
into a new vertex, is the inverse of the previous one. Because not
all cycles can be so contracted, stronger hypotheses must be made
on the graph $G$ and on the cycles considered. To state those conditions,
we first need some definitions. 
\begin{defn}
A \textbf{$k$-regular} graph is a graph where all vertices have the
same degree $k$. A \textbf{quasi $k$-regular} graph is a graph where
at most one vertex has a degree different than $k$.\end{defn}
\begin{lem}
An exactly $k$-edge connected graph $G$ which has only trivial cuts
is quasi $k$-regular.\end{lem}
\begin{proof}
Suppose there exists two vertices $u$ and $v$ of degree greater
than $k$. There exists a cut separating $u$ and $v$ and this cut
cannot be trivial.\end{proof}
\begin{prop}
[Cycle contraction]\label{lem:Cycle-Contraction}Let $G=(V,E)$ be
a quasi 3-regular exactly 3-edge-connected and biconnected graph.
Let $C$ be a chordless cycle of any length which contains the vertex
of higher degree. The graph $G'$ is contructed from $G$ by collapsing
the cycle $C$ into one vertex. Formaly, if $2\le d'\le d$, $u_{1},\ldots,u_{d'}$
are the vertices of $C$, and $v_{1},\ldots,v_{d}$ are the vertices
adjacent to $C$ in $G$, we construct $G'=(V',E')$ with\begin{align*}
V'= & V\setminus\{u_{1},\ldots,u_{d'}\}\cup\{u\}\\
E'= & E\cup\{(u,v_{i}):1\leq i\leq d\}\setminus\left(\left\{ (u_{i,}v_{i}):1\leq i\leq d'-1\right\} \cup\left\{ (u_{d'},v_{i}):d'\leq i\leq d\right\} \right)\,.\end{align*}
Then $G'$ is exactly 3-edge connected and quasi 3-regular.\end{prop}
\begin{proof}
In $G'$, the degree of every vertex in $V\setminus\{u_{1},\ldots,u_{d'}\}$
is the same as its degree in $G$. Given that only $u_{d'}$ in $G$
may have a degree different than 3, $G'$ is also quasi 3-regular.
In particular, there cannot be more than 3 edge-disjoint paths between
any pair of vertices in $G'$. For any pair of vertices $x,y$ not
on $C$, any $x-y$ path in $G$ can be changed into a path in $G'$
where the eventual edges on the cycle are removed. Furthermore, edge-disjoint
paths in $G$ are still edge-disjoint in $G'$. So there are 3 edge-disjoint
$x-y$ paths in $G'$. For any $x$ not on $C$, any $x-u_{1}$ path
in $G$ can be changed into a $x-u$ path in $G'$. So there are 3
edge-disjoint paths in $G'$ between any pair of vertices.
\end{proof}

\subsection{Main Result: Synthesis\label{sub:Synthesis}}

\begin{figure}
\begin{centering}
\includegraphics[clip,width=0.25\textwidth]{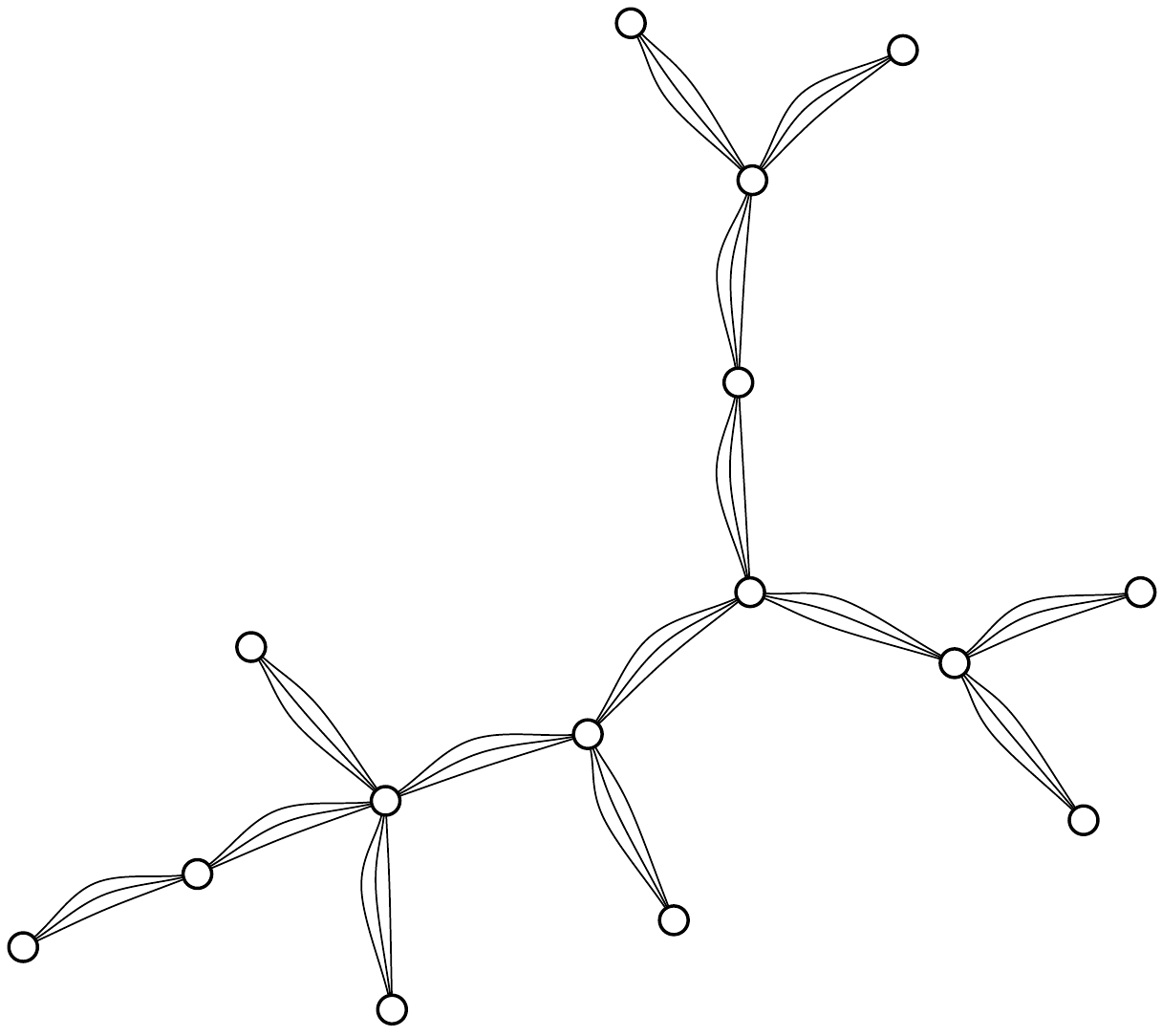}
\par\end{centering}

\caption{\label{fig:A-3-thick-tree.}A 3-thick tree.}

\end{figure}

We now proceed to the synthesis of exactly 3-edge-connected graphs.
The proof will rely on results from Section \ref{sub:Existence-of-Non-articulation-Cycles},
where the main technical arguments are given. 
\begin{defn}
The \textbf{dumbbell graph} consists of 2 vertices and 3 parallel
edges between them. It is exactly 3-edge-connected.\end{defn}
\begin{thm}
[Exactly 3-edge-connected synthesis]\label{thm:synthesis}Any exactly
3-edge-connected graph $G$ is obtained from dumbbell graphs and the
following operations: cycle expansion and block gluing.\end{thm}
\begin{proof}
We proceed by induction on the order $r$ of $G$. The only exactly
3-edge-connected graph of order 2 is the dumbbell and the induction
hypothesis holds for $r=2$. Suppose the theorem holds for any graph
of order $j\le r$. Let $G$ be an exactly 3-edge connected graph
of order $r+1$. We apply the following tests:
\begin{enumerate}
\item If $G$ has an articulation point, by Corollary~\ref{cor:Block-splitting},
each block is exactly 3-edge-connected and of order less than $r+1$.
Apply the induction to each block; $G$ is obtained by block gluing
of its blocks.
\item $G$ is biconnected and of order at least 3, so by Theorem~\ref{thm:Collapsable-cycles}
it contains a quasi 3-regular chordless cycle $C$ which, when collapsed,
does not create an articulation point and preserves exact 3-edge-connectivity.
Let $G'$ be obtained from $G$ by collapsing $C$. The order of $G'$
is less than $r+1$. Apply the induction to $G'$; $G$ is obtained
by cycle expansion of $G'$.
\end{enumerate}
\end{proof}
\begin{cor}
[3-thick tree structure]Every exactly 3-edge-connected graph can
be obtained from a 3-thick tree (Figure~\ref{fig:A-3-thick-tree.})
and block-respecting cycle expansions.\end{cor}
\begin{proof}
The block gluing and cycle expansion operations are independent of
each other and can be applied in any order. First apply a sequence
of block gluing of dumbbells to duplicate the block structure of the
graph and get a 3-thick tree. Then apply block-respecting cycle expansion
operations.
\end{proof}

The main difficulty in the proof of Theorem~\ref{thm:synthesis}
is that collapsing an arbitrary cycle could create an articulation
point, but cycle expansion (Proposition~\ref{lem:Cycle-expansion})
and cycle contraction (Proposition~\ref{lem:Cycle-Contraction})
are inverse operations of each other only as long as the contraction
does not create an articulation point. This difficulty will be dealt
with in the next section.

\subsection{Existence of Non-articulation Cycles\label{sub:Existence-of-Non-articulation-Cycles}}

In this section, we shall prove the existence in an exactly 3-edge-connected
graph of a quasi 3-regular chordless cycle which, when collapsed,
does not create an articulation point and preserves exact 3-edge-connectivity.
This will culminate in Theorem~\ref{thm:Collapsable-cycles}, which
was used to prove that all exactly 3-edge-connected graphs can be
constructed from block-gluing and cycle expansion (Theorem~\ref{thm:synthesis}). 

Several previous works \cite{Thomassen81,Thomassen81a,Egawa02} deal
with the problem of finding non-separating induced cycles in a graph,
i.e.\@ a chordless cycle $C$ in $G$ such that $G\setminus C$ is
connected. Thomassen's results \cite{Thomassen81} for 2-connected
graphs of minimum degree 3 would guarantee the existence of such a
cycle. But collapsing such cycle does not necessarily preserve exact
3-edge-connectivity as required for the synthesis, and stronger conditions
must be enforced. We introduce some definitions and operations relating
to articulation points in exactly 3-edge-connected graphs. 
\begin{defn}
Let $H$ be a subgraph of a connected graph $G=(V,E)$ and let $V_{1},\ldots,V_{k}$
be the sets of vertices of the $k$ connected components of $G\setminus H$.
The partition $V=H\cup V_{1}\cup\ldots\cup V_{k}$ induced by $H$
is the \textbf{$H$-partition} and its size is $k$ ($k\ge0$).
\end{defn}

\begin{figure}
\begin{centering}
\includegraphics[clip,width=0.7\textwidth]{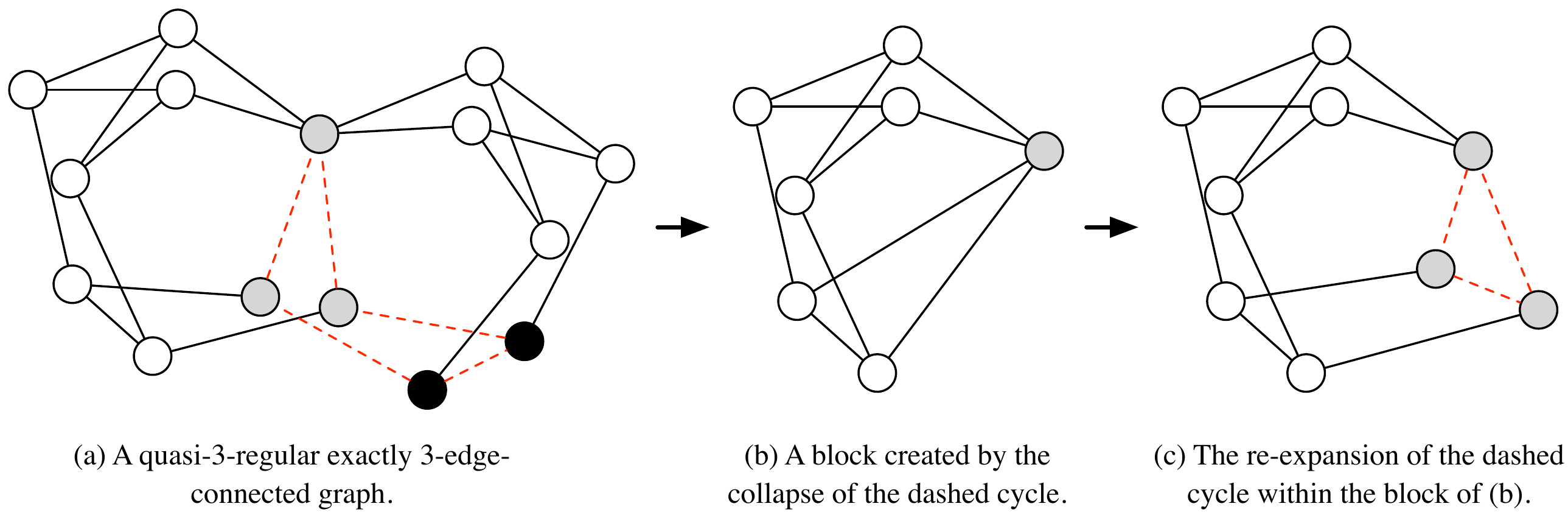}
\par\end{centering}

\caption{\label{fig:contra-exp}An example contraction-expansion of a cycle.}

\end{figure}

\begin{defn}
An \textbf{articulation cycle} in $G=(V,E)$ is a chordless cycle
$C$ with a $C$-partition of size greater than 1.
\end{defn}
An equivalent definition is that $C$ is an articulation cycle of
a biconnected graph $G$ if $G'$ obtained from $G$ by contracting
$C$ into a single vertex has more than 1 biconnected component.
\begin{defn}
Let $G$ be a graph, $C$ be a chordless cycle in $G$ and $V$ be
a block of the $C$-partition. The \textbf{contraction-expansion}
of $C$ with respect to $V$ is a graph $G'$ obtained from $G$ by
the following operations, applied sequentially:
\begin{itemize}
\item let $G'$ be the graph induced by $V\cup C$
\item smooth out the vertices of degree two in $G'$
\end{itemize}
\end{defn}
The vertices of degree two that are smoothed out (replaced by edges)
are the attachment vertices of the blocks other than $V$ in the $C$-partition.
The name contraction-expansion is justified by the proof of Lemma~\ref{lem:contraction-expansion-properties}.
See Figure~\ref{fig:contra-exp} for an example of such a contraction-expansion.

\begin{lem}
\label{lem:contraction-expansion-properties}Let $G$ be an exactly
3-edge-connected graph, $C$ be a chordless cycle in $G$ and $V$
be a block of the $C$-partition. Then $G'$ the contraction-expansion
of $C$ with respect to $V$ is exactly 3-edge-connected. Furthermore,
if $G$ is 3-regular then so is $G'$. If $G$ is quasi 3-regular
and $C$ contains the high-degree vertex, then $G'$ is quasi 3-regular
and $C'$ contains the high-degree vertex of $G'$.\end{lem}
\begin{proof}
$G'$ can also be constructed by the following operations: contract
the cycle $C$, let $G'$ be the block which contains $V$ and expand
in $G'$ the contracted cycle into a cycle $C'$. The first part of
the theorem statement follows from the properties of the cycle contraction
(Proposition~\ref{lem:Cycle-Contraction}), block gluing (Corollary~\ref{cor:Block-splitting})
and cycle expansion (Proposition~\ref{lem:Cycle-expansion}) operations.
After a cycle expansion, all the vertices but one in the newly created
cycle have degree 3. By construction, a vertex in $C'$ has a degree
less or equal to the corresponding vertex in $C$. Hence, the 3-regularity
and quasi 3-regularity properties are preserved\end{proof}
\begin{defn}
A \textbf{collapsible} cycle is a quasi 3-regular, non-articulation,
chordless cycle which preserves exact 3 edge-connectivity when collapsed.\end{defn}
\begin{lem}
\label{lem:regular-imply-simple}Let $G$ be a 3-edge-connected graph.
If $G$ is 3-regular, then it is a simple graph. If $G$ is quasi
3-regular, the high degree vertex $h$ is an end point of any double
edge.\end{lem}
\begin{proof}
If $G$ is 3-regular and has a double edge $(u,v)^{2}$ then $u$
(resp.\@ $v$) has only one other neighbor beside $v$ (resp.\@
$u$), say $u'$ (resp.\@ $v'$). The set $S=\left\{ \left(u,u'\right),\left(v,v'\right)\right\} $
is exactly the edge cut $\left\langle \left\{ u,v\right\} ,G\setminus\left\{ u,v\right\} \right\rangle $
and has size 2, which is a contradiction. So $G$ must be a simple
graph. If $G$ is quasi 3-regular and $\left(u,v\right)^{2}$ is a
double edge not incident to $h$, then the same contradiction as above
arises. So any double edge must be incident to $h$.\end{proof}
\begin{lem}
\label{lem:Cycle-wo-(v,w)-w-u}Let $G$ be a 3-regular, 3-edge-connected
graph, $(v,w)$ be an edge in $G$ and $u$ be a vertex in $G$ distinct
from $v$ and $w$. Then $G$ contains a chordless cycle that contains
edge $(v,w)$ and does not contain $u$.\end{lem}
\begin{proof}
Because there are 3 edge-disjoint paths between $v$ and $w$, edge
$(v,w)$ is on at least 2 cycles $C_{1}$ and $C_{2}$. The cycles
are vertex disjoint except for $v$ and $w$ because $G$ is 3-regular.
Hence, one of $C_{1}$ and $C_{2}$ does not contain $u$. If it contains
chords, it can be short-circuited to change it to a chordless cycle,
which still does not contain $u$.\end{proof}
\begin{lem}
\label{lem:Quasi-regular-cycle-wo-(v,w)-w-u}Let $G$ be a quasi 3-regular,
3-edge-connected graph, $\left(v,w\right)$ be an edge in $G$ where
$v$ or $w$ is the high degree vertex and $u$ be a vertex in $G$
distinct from $v$ and $w$. Suppose also that $G$ has at most one
double edge and $u$ is an end point of the double edge if any. Then
$G$ contains a chordless cycle that contains edge $\left(v,w\right)$
and does not contain $u$.\end{lem}
\begin{proof}
By an identical proof as in Lemma~\ref{lem:Cycle-wo-(v,w)-w-u},
there is a cycle $C$ which contains $\left(v,w\right)$ and not $u$.
This cycle cannot have any double edge, and if it contains any chords
it can be short-circuited.\end{proof}
\begin{lem}
\label{lem:4-bicolored-edges}Let $G$ be a quasi 3-regular, 3-edge-connected
graph and let $C$ be a chordless articulation cycle with a partition
of size $k>1$ which contains the high degree vertex $h$. Each vertex
of $C$ connects to exactly one block. Color those vertices that connect
to a given block $V$ red and color the others blue. Then there are
at least 4 edges on $C$ that connect a red vertex to a blue vertex.\end{lem}
\begin{proof}
Because C is a cycle, there must be an even number of red-blue edges.
If there were no such edges, then $k$ would be 1. If there were only
2 such edges, they would form a 2-edge-cut. Hence $k\ge4$.
\end{proof}
The following lemmas and subsequent theorem contain the main technical
arguements required to ensure that a reversable cycle contraction
can be performed on any biconnected, quasi 3-regular, exactly 3-connected
graph. Theorem~\ref{thm:Collapsable-quasi-regular}, which uses Lemma~\ref{lem:Collapsable-cycles-3-regular},
will handle one of the cases of the induction proof of Theorem~\ref{thm:synthesis}.
In these proofs, we often use the the following technique: given an
exactly 3-edge-connected graph $G$, we create a new exactly 3-edge-connected
graph $G'$ via contraction-expansion of some cycle $C$ in $G$.
We say a cycle $Z\ne C$ in $G$ is the \textbf{corresponding cycle}
to a cycle $Z'$ in $G'$ if $Z$ uses every edge of $Z'$ that is
present in $G$, and if $Z'$ uses an edge $(u,v)$ that is not in
$G$ then $u,v$ are on cycle $C$ and $Z$ uses a path on $C$ between
$u$ and $v$.

\begin{lem}
[Collapsible cycles in 3-regular graphs]\label{lem:Collapsable-cycles-3-regular}Let
$G=(V,E)$ be a biconnected exactly 3-edge connected 3-regular graph,
let $(v,w)\in E$ be an edge and $u\in V$ be a vertex distinct from
$v$ and $w$. Then $G$ contains a collapsible cycle which contains
edge $(v,w)$ and does not contain vertex $u$.\end{lem}
\begin{proof}
By Lemma~\ref{lem:regular-imply-simple}, $G$ is a simple graph.
We proceed by induction on $n$, the number of vertices of $G$. It
is true when $n=4$ ($K_{4}$, the complete graph on 4 vertices, is
the only such graph).

Let $C$ be a chordless cycle in $G$ that contains $(v,w)$ and doesn't
contain $u$ and that, among those, minimizes $k$ the size of the
partition. The existence of such a cycle is guaranteed by Lemma~\ref{lem:Cycle-wo-(v,w)-w-u}.
If $k=1$, then $C$ is the desired cycle.

For purpose of contradiction, suppose that $k>1$. Vertex $u$ is
in only one block. Let $V$ be a block in the $C$-partition which
does not contain $u$. Color in $C$ the $V$-adjacent vertices red
and the non-$V$-adjacent vertices blue. Because $G$ is 3-regular
and 3-edge-connected, $C$ must contain at least 3 red vertices and
at least 3 blue vertices. By Lemma~\ref{lem:4-bicolored-edges},
there are at least 4 red-blue edges. Therefore, there is a path $P$
in $C$ that starts and ends at red vertices and uses edge $(v,w)$
but does not use some red-blue edge $e$. Let $x$ and $y$ be respectively
the red vertex and blue vertex of edge $e$, with associated blocks
$V$ and $W_{y}$.

Let $G'$ be the contraction-expansion of $C$ with respect to $V$.
Because vertex $x$ was colored red it is still present in $C'$ and
to the path $P$ in $C$ corresponds an edge $P'$ of $C'$. By Lemma~\ref{lem:contraction-expansion-properties},
$G'$ is 3-regular and exactly 3-edge-connected. $G'$ has fewer vertices
than $G$ and by induction there exists a collapsible cycle $Z'$
which contains the edge $P'$ and does not contain vertex $x$. Because
$x\in C'$, the cycles $Z'$ and $C'$ are distinct. Because $P'\subset Z'\cap C'$,
$Z'$ and $C'$ are not disjoint.

Let $Z$ be the cycle in $G$ corresponding to $Z'$ and let $k'$
be the size of the $Z$-partition. Because $Z$ is contained in $C\cup V$,
if $a$ and $b$ are two vertices in a block $W$ (distinct from $V$)
of the $C$-partition, $a$ and $b$ must belong to the same block
$W'$ in the $Z$-partition. Moreover, $Z'$ was not an articulation
cycle in $G'$, hence if $a,b$ are two vertices in $V\setminus Z$,
they must belong to the same block of the $Z$-partition. Hence, there
is a surjective function from the $C$-partition to the $Z$-partition
and $k'\le k$. If $W$ is a block in the $C$-partition, let note
$W'$ the unique block in the $Z$-partition containing $W$.

By construction, vertex $x$ is not contained in the cycle $Z$, so
the red-blue edge $e$ is also not contained in the cycle $Z$. Also,
in $G\setminus Z$ there is a path between $V'$ (the block containing
$x$) and $W'_{y}$ (the block containing $y$). In other words, $V'=W'_{y}$,
the function from the $C$-partition to the $Z$-partition is not
injective and $k'<k$. This contradicts the minimality of $k$ and
concludes the induction step. \end{proof}
\begin{thm}
[Collapsible cycles in quasi 3-regular graphs]\label{thm:Collapsable-quasi-regular}Let
$G=(V,E)$ be a biconnected, exactly 3-edge-connected, quasi 3-regular
graph with high-degree vertex $h$ and let $u$ be a vertex of $G$
distinct from $h$. Then $G$ contains a collapsible cycle which contains
$h$ and does not contain $u$.\end{thm}
\begin{proof}
If $G$ has a double edge which is not incident to $u$ (in particular
if there are two or more double edges), then this double edge is a
collapsible cycle, which contains $h$ by Lemma~\ref{lem:regular-imply-simple},
and avoids $u$. So we assume that $G$ has at most one double edge,
and if it does, $u$ is an end point of the double edge. We proceed
similarly to Lemma~\ref{lem:Collapsable-cycles-3-regular} and prove
the statement by induction on $n$, the number of vertices of $G$.
The theorem is true when $n=4$ (again, $K_{4}$ is the only such
graph).

Let $C$ be a chordless cycle in $G$ that contains $h$ and does
not contain $u$ and that, among those, minimizes the size $k$ of
the partition. The existence of such a cycle is guaranteed by Lemma~\ref{lem:Quasi-regular-cycle-wo-(v,w)-w-u}.
If $k=1$, then $C$ is the desired cycle. For purpose of contradiction,
suppose that $k>1$. Let $V$ be a block in the $C$-partition which
does not contain $u$. Color in $C$ the $V$-adjacent vertices red
and the non-$V$-adjacent vertices blue. By Lemma~\ref{lem:4-bicolored-edges},
there are at least 4 red-blue edges. The two cases, where $h$ is
adjacent to $V$ or not, are now discussed.

If $h$ is adjacent to $V$, then $h$ is colored red. Let $e$ be
any red-blue edge which does not have $h$ as one of its end point.
Let $x$ and $y$ be respectively the red vertex and blue vertex of
edge $e$, with associated blocks $V$ and $W_{y}$. Consider $G'$
the contraction-expansion of $C$ with respect to $V$. $G'$ contains
$h$ as $h$ is red and, by Lemma~\ref{lem:contraction-expansion-properties},
$G'$ is quasi 3-regular and exactly 3-edge-connected. $G'$ has fewer
vertices than $G$ and by induction there exists a collapsible cycle
$Z'$ which contains $h$ and does not contain $x$. The corresponding
cycle $Z$ in $G$ contains $h$ and does not contain edge $e$.

If $h$ is not adjacent to $V$ then $h$ is colored blue. Let $P$
be a $v-w$ path in $C$ between two red vertices such that $h$ is
an internal of $P$ and all the internal vertices of $P$ are blue.
Let $e$ be any red-blue edge which does not have $v$ or $w$ as
an end point and let $x$ and $y$ be respectively the red vertex
and blue vertex of edge $e$, with associated blocks $V$ and $W_{y}$.
Let $G'$ be the contraction-expansion of $C$ with respect to $V$.
To the path $P$ now corresponds the edge $(v,w)$ in $C'$ and $x$
is in $G'$ because it is red. The high degree vertex $h$ of $G$
is not in $G'$ because it is blue, hence, by Lemma~\ref{lem:contraction-expansion-properties},
$G'$ is 3-regular and exactly 3-edge-connected. By Lemma~\ref{lem:Collapsable-cycles-3-regular},
there exists a cycle $Z'$ in $G'$ which contains the edge $(v,w)$
and does not contain vertex $x$. The corresponding cycle $Z$ in
$G$ contains $P$ and does not contain $x$. Hence cycle $Z$ contains
$h$ and does not contain edge $e$.

In both cases, a cycle $Z$ containing $h$ and not containing a red-blue
edge $e$ exists. An identical argument as in Lemma~\ref{lem:Collapsable-cycles-3-regular}
shows that the $Z$-partition has size $k'<k$, which contradicts
the minimality of $k$.
\end{proof}
We now prove the existence of a collapsible cycle in the general case,
where the only conditions on $G$ are that it is biconnected and exactly
3-edge-connected. We start with a lemma that guarentees that such
a graph has vertices of degree 3. The trees, which are the exactly
1-edge-connected graph, have two leaves (vertices of degree 1), and
the exactly 2-edge-connected graphs also have two vertices of degree
2. Similarly, we have the following property for an $k\ge3$.
\begin{lem}
\label{lem:Min-degree-3}Let $G$ be a biconnected, exactly $k$-edge-connected
graph of order $\geq k$. Then $G$ has at least 2 vertices of degree
$k$.\end{lem}
\begin{proof}
We proceed by induction on the number of non-trivial minimum cuts
in $G$. If $G$ has no non-trivial minimum cuts, then it is quasi
$k$-regular. If $G$ is $k$-regular, it has at least 2 vertices
of degree $k$. If $G$ is quasi $k$-regular with a high degree vertex,
it has at least 3 vertices and therefore has at least 2 vertices of
degree $k$.

Let $S=\langle V_{1},V_{2}\rangle$ be a non-trivial minimum cut and
let $G_{1}$ and $G_{2}$ be the vertex splitting graphs induced by
$S$. Call $x_{1}$ and $x_{2}$ the new vertices ($V(G_{i})=V_{i}\cup\{x_{i}\}$).
Suppose $T$ were a non-trivial minimum cut in $G_{1}$. Construct
a corresponding non-trivial minimum cut $T'$ in $G$ by changing
any edge used by $T$ that is adjacent to $x_{1}$ to the corresponding
edge in $S$. So to any non-trivial minimum cut of $G_{1}$ or $G_{2}$
corresponds a distinct non-trivial minimum cut in $G$. But no non-trivial
cut in $G_{1}$ or $G_{2}$ corresponds to the cut $S$ (it would
be a trivial cut in $G_{1}$ and $G_{2}$ ). So both $G_{1}$ and
$G_{2}$ have fewer non-trivial minimum cuts than $G$. By induction
$G_{1}$ and $G_{2}$ have 2 vertices of degree $k$, including $x_{1}$
and $x_{2}$. So $G$ is obtained as a vertex gluing of $G_{1}$ and
$G_{2}$ and has 2 vertices of degree $k$.\end{proof}
\begin{thm}
[Collapsible cycles]\label{thm:Collapsable-cycles}Let $G$ be a
biconnected, exactly 3-edge-connected graph of order at least 3 and
$u$ a vertex of degree 3. Then $G$ contains a collapsible cycle
which does not contain $u$.\end{thm}
\begin{proof}
By induction on the number of non-trivial minimum cuts in $G$. If
$G$ has no non-trivial minimum cuts, then it is quasi 3-regular.
By Theorem~\ref{thm:Collapsable-quasi-regular} it contains a collapsible
cycle which avoids $u$.

Assume $G$ has non-trivial minimum cuts. Let $S=\langle V_{1},V_{2}\rangle$
be a non-trivial minimum cut, $\left|S\right|=3$, and let $G_{1}$
and $G_{2}$ be the vertex splitting graphs induced by $S$. Assume
without loss of generality that $G_{1}$ does not contain $u$. By
the same argument as in Lemma~\ref{lem:Min-degree-3}, $G_{1}$ has
less non-trivial minimum cuts than $G$ and by induction, $G_{1}$
has a collapsible cycle $C_{1}$ which avoids $x_{1}$.

$C_{1}$ in $G$ is not an articulation cycle. Moreover, contracting
$C_{1}$ in $G$ also preserves exact 3-edge connectivity: collapsing
a cycle does not destroy any path between the remaining vertices,
it does not create 4 paths in $V_{1}$ by construction, and it does
not create 4 paths between $V_{1}$ and $V_{2}$ because of the minimum
cut $S$. Therefore, $C_{1}$ is a collapsible cycle of $G$ which
does not contain $u$.
\end{proof}

\section{Other properties of exactly $k$-edge-connected graphs}

\subsection{Number of Operations}

Let $G=(V,E)$ be an exactly 3-edge-connected graph with $n_{G}=\left|V\right|$
and $m_{G}=\left|E\right|$, and let $B_{G}$ be the number of blocks
in $G$. The number of synthesis operations needed to generate $G$
is determined by $m_{G},n_{G}$ and $B_{G}$.
\begin{prop}
If $G=(V,E)$ is an exactly 3-edge-connected graph, then \begin{equation}
m_{G}=n_{G}+2B_{G}+E_{G}-1\,,\label{eq:nb-cycle-expansions}\end{equation}
where $E_{G}$ is the number of cycle expansions in a synthesis of
$G$.\end{prop}
\begin{proof}
A dumbbel graph has 3 edges, 2 vertices and one block. It satisfies
$m=n+B$. Let $T$ be 3-thick tree satisfying $m_{T}=n{}_{T}+2B_{T}-1$,
and let $T'$ the block addiction of $T$ with a dumbbell. We have
$m_{T'}=m_{T}+3$, $B_{T'}=B_{T}+1$ and $n_{T'}=n_{T}$. So $T'$
satisfies $m_{T'}=n{}_{T'}+B_{T'}$. The cycle expansion operation
on $G$ with a cycle of $d$ vertices adds $d-1$ vertices and $d$
edges. So if $G$ is obtained from a 3-thick tree by $E_{G}$ cycle
expansions, relation \ref{eq:nb-cycle-expansions} holds. To create
$B_{G}$ blocks, $B_{G}-1$ block gluing operations are needed. So
the the number of operations in a synthesis of $G$ is $N_{G}=B_{G}-1+E_{G}=m_{G}-n_{G}-B_{G}$.
\end{proof}
%
{}

\subsection{Minimum Exactly Connected Graphs}

A natural requirement for network design is to use the fewest edges
possible. A $k$-vertex-connected or $k$-edge-connected graph with
$n$ vertices has at least $\left\lceil \frac{kn}{2}\right\rceil $
edges. We say a graph is \textbf{minimum} if it has exactly that many
edges. The Harary~\cite{Gross2006} graph $H_{k,n}$ is a special
graph which has $n$ vertices, is $k$-vertex-connected and $k$-edge-connected,
and is minimum. They are, in addition, quasi $k$-regular, which implies
that they are also exactly $k$-edge-connected. In fact, we have following
relationships between minimum, almost $k$-regular, and exactly $k$-edge-connected
graphs:
\begin{defn}
A graph $G$ is \textbf{almost} $k$-regular if it is quasi $k$-regular
and has maximum degree $\leq k+1$.\end{defn}
\begin{prop}
Let $G$ be a graph. The following assertions are equivalent:
\begin{enumerate}
\item $G$ is $k$-edge-connected and minimum.
\item $G$ is $k$-edge-connected and almost $k$-regular.
\item $G$ is exactly $k$-edge-connected and almost $k$-regular.
\end{enumerate}
\end{prop}
We therefore have that, among the $k$-edge-connected graphs, the
minimum graphs are a strict subset of the exactly connected graphs,
which in turn are a strict subset of the edge-minimal graphs. The
minimum, exactly 3-edge-connnected graphs can be generated by disallowing
a cycle expansion if it would lead to $\sum_{v}\textrm{deg}(v)>3n+1$.

\subsection{Planar graphs}

Planar, exactly 3-edge-connected graphs can be generated with a slightly
modified synthesis. If $G$ is a planar, exactly 3-edge-connected
graph and $u$ is a vertex, a drawing of $G$ on the plane defines
a natural order on the neighbors of $u$ (for example, clockwise traversal
of the edges incident to $u$). An \textbf{order-preserving cycle
expansion} is a cycle expansion of $u$ where the vertices in the
cycle are attached to the neighbors of $u$ in the natural order.
$G'$ obtained from $G$ by an order-preserving cycle expansion is
also a planar graph. Conversely, if $G$ is planar and $C$ is a collapsible
cycle, $G'$ obtained by collapsing $C$ is also planar.
\begin{thm}
[Planar exactly 3-edge-connected synthesis]A graph is planar and
exactly 3-edge-connected if and only if it can be obtained from dumbbell
graphs and the folowing operations: order-preserving cycle expansion
and block gluing.
\end{thm}

\section{Conclusion}

We have given an operation (cycle expansion) by which every exactly
3-edge-connected, biconnected multigraph can be generated. The synthesis
allows for the generation of all fair, 3-edge-connected networks,
and preserves the block tree structure underlying such graphs. A natural
question is whether this synthesis operation can be used to design
fair, robustly connected networks that satisfy other desired properties.

\bibliographystyle{amsplain}
\bibliography{bib}

\end{document}